\documentclass[11pt]{amsart}

\usepackage[english]{babel}
\usepackage{amssymb}
\usepackage{fourier}
\usepackage{mathrsfs}
\usepackage{enumerate}
\usepackage{color}
\usepackage{ifpdf}
\usepackage{tikz}
\usepackage{tikz-cd}
\usetikzlibrary{decorations.pathmorphing,arrows, intersections}
\usepackage[initials]{amsrefs}
\usepackage{enumitem}

\usepackage{caption}
\usepackage{subcaption}
\usepackage{comment}

\newcommand{\bfH}{{\mathbf H}}

\newcommand{\SL}{\operatorname{SL}}
\newcommand{\PSL}{\operatorname{PSL}}

\newtheorem{theorem}{Theorem}[section]
\newtheorem{thm}[theorem]{Theorem}

\newtheorem{lemma}[theorem]{Lemma}

\newtheorem{defn}[theorem]{Definition}

\newtheorem{example}[theorem]{Example}

\newtheorem{remark}[theorem]{Remark}

\newtheorem{question}[theorem]{Question}

\numberwithin{equation}{section}

\begin{document}
\title[Marked length varieties and arithmeticity]{Marked length varieties and arithmeticity}

\author{Yanlong Hao}
\address{University of Michigan, Ann Arbor, Michigan, USA}
\email{ylhao@umich.edu}
 
\date{\today}
    
\maketitle

\begin{abstract}
    This paper gives examples of non-arithmetic surfaces with marked length variety rigidity.
\end{abstract} 
\section{Introduction}
Let $(M,g)$ be a closed negatively curved Riemannian manifold. In this setting, every free homotopy class of closed curves contains a unique closed geodesic. The \textbf{marked length spectrum} map is the function that assigns to each free homotopy class the length of its corresponding closed geodesic, namely
\[\ell_g:\pi_1(M)\to \mathbb{R}.
\] For non-compact negatively curved manifolds, the marked length spectrum can be extended by assigning length zero to parabolic elements.

The marked length spectrum rigidity problem, proposed in \cite{burns1985manifolds}*{Conjecture 3.1}, asks whether the marked length spectrum of a closed negatively curved manifold uniquely determines its isometry class. This conjecture has been verified in several cases: for surfaces by Otal and independently by Croke \cites{otal1990spectre, croke1990rigidity}, and for closed locally symmetric manifolds by Hamenst{\"a}dt \cite{hamenstadt1997cocycles}, using the entropy rigidity theorem of Besson-Courtois-Gallot \cite{besson1995entropies}. 

For general negatively curved metrics, Katok's work \cite{katok1988four} gives an affirmative answer to the marked length spectrum rigidity problem within a fixed conformal class. Moreover, a recent result of Guillarmou and Lefeuve \cite{guillarmou2019marked} shows that the conjecture holds locally. 

Related results for the conjecture in non-manifold settings include flat surface metrics \cite{bankovic2018marked}, one-dimensional spaces \cite{constantine2019marked}, CAT(0) cube complexes \cite{beyrer2021cross}.

In \cite{hao2022marked}, we investigated a related aspect of this problem. We define $P_g$ to be the set of pairs of elements in $\Gamma=\pi_1(M)$ whose corresponding closed geodesics have equal length; this set is
called the \textbf{marked length pattern} of $(M,g)$.
One of the main results in \cite{hao2022marked} establishes marked length pattern rigidity for closed arithmetic manifolds and for non-compact arithmetic surfaces:

\textbf{Let $(M, g_0)$ be a closed arithmetic Riemannian manifold, and let $g$ be an arbitrary negatively curved Riemannian metric on $M$. If $P_{g_0}\subset P_{g}$, then $g=\lambda g_0$ for some $\lambda>0$.}

This motivates the study of marked length patterns. In particular, one is led to the following question:
\begin{question}\label{inverse}
If a closed negatively curved manifold $(M, g)$ exhibits marked length pattern rigidity, must $M$ be arithmetic?  
\end{question}

In this paper, we introduce a stronger variant of the marked length pattern, called the \textbf{marked length variety} (see Definition~\ref{def: MLV}). We give examples of non-arithmetic surfaces exhibiting marked length variety rigidity. Nevertheless, Question~\ref{inverse} remains open.

Before presenting our results, we fix some notation and recall several definitions.

Let $(M, g)$ be a complete, pinched, negatively curved manifold of finite volume, and let $\ell_g$  denote its marked length spectrum. For a closed geodesic $\gamma$, we define its \textbf{trace} by \[\mathrm{Tr}_g(\gamma)=2 \cosh\left(\frac{\ell_g(\gamma)}{2}\right).\] We further define the \textbf{extended fundamental group} of $M$ by \[\bar{\pi}_1(M)=\{\pm1\}\times\pi_1(M)\] as a set, and the \textbf{signed trace} on $\bar{\pi}_1(M)$ by $
\mathrm{Tr}_g(\pm 1, [\gamma])=\pm\mathrm{Tr}_g(\gamma)$. Here, a conjugacy class $[\gamma]$ of $\pi_1(M)$ is identified with the closed geodesic that represents it.

When $(M, g)$ is a hyperbolic surface, the fundamental group $\pi_1(M)$ is a subgroup of $\PSL_2(\mathbb{R})$, and $\bar\pi_1(M)$ can be identified with the preimage of $\pi_1(M)$ under the standard projection \[p: \SL(2,\mathbb{R})\to\PSL(2,\mathbb{R}).\] In this case, $\mathrm{Tr}_g$ coincides with the usual matrix trace, and the group structure on $\bar{\pi}_1(M)$ is induced from that of $p^{-1}(\pi_1(M))$.

Let $F\in \mathbb{Q}[X_1,X_2,\ldots,X_n]$ be a polynomial.

\begin{defn}\label{def: MLV}
The marked length variety of $F$ associated with $g$, denoted $V_F^g$, is the set of tuples $(\gamma_1,\gamma_2,\cdots,\gamma_n)\in \bar{\pi}_1(M)^n$ such that 
$$F(\mathrm{Tr}_g(\gamma_1),\mathrm{Tr}_g(\gamma_2),\cdots,\mathrm{Tr}_g(\gamma_n))=0.$$
\end{defn}

\begin{example}
For $H=X_1-X_2$, the variety $V_H^g$ coincides with the marked length pattern associated with $g$.
\end{example}

Given a set $S$ of polynomials, define $V_S^g=\prod_{F\in S}V_F^g$.
\begin{defn}
    A manifold $(M,g)$ is said to be marked length variety rigid if there exists a finite set $S$ of polynomials such that for any complete, pinched, negatively curved metric $g'$ on $M$, if $V_S^g\subset V_S^{g'}$, then $g'=\lambda g$ for some $\lambda>0$. 
\end{defn}

As an illustration, we provide an example of a non-arithmetic surface exhibiting marked length variety rigidity, where $S$ consists of two polynomials.

Let $S_{1,1}$ be a surface of genus one with a single puncture, and let $a$, $b$ be standard generators of $\pi_1(S_{1,1})$. There exists a unique hyperbolic structure $g$ on $S_{1,1}$ satisfying
\[
\mathrm{Tr}_g(a)=\mathrm{Tr}_g(b),\quad \mathrm{Tr}_g(a^2)=\mathrm{Tr}_g(a^2b).
\]

In section~\ref{Proofs}, we show that this metric $g$ is marked length variety rigid with respect to the set of polynomials $\{X_1-X_2,\  X_1X_2-X_3-X_4\}$. We also construct an example of a marked length variety rigid closed surface in the same section.

\medskip
We also present a method for detecting marked length variety rigid surfaces, beginning with a definition.

\begin{defn}
     A \textbf{Salem number} $x$ is a real algebraic integer satisfying $x>1$ such that all of its Galois conjugates except $\frac{1}{x}$ lie on the unit circle.

    A real algebraic number $x$ is called a \textbf{geometric Salem number} if $x>2$ and all of its Galois conjugates are real with absolute values less than 2. 
\end{defn}

\begin{remark}
    Let $x$, $y\in \mathbb{R}$ satisfy $x>2$, $y>1$, and $y+\frac{1}{y}=x$. Then $x$ is a geometric Salem number if and only if $y$ is a Salem number. This follows from the following observations: for any Galois conjugate $\sigma(y)$ of $y$, 
    \begin{enumerate}
        \item $\sigma(y)$ is real if and only if $\sigma(x)$ is a real and $|\sigma(x)|\geq 2$. 
        \item $\sigma(y)$ lies on the unit circle and is not real if and only if $\sigma(x)$ is real and $|\sigma(x)|<2$.
        \item $\sigma(x)=x$ if and only if $\sigma(y)=y$ or $\sigma(y)=\frac{1}{y}$.
    \end{enumerate}
\end{remark}
\begin{thm}\label{A}
   Let $(\Sigma,g)$ be an orientable, complete, finite-volume hyperbolic surface, with $\Gamma=\pi_1(\Sigma)$ a Fuchsian lattice. Suppose $\Gamma$ admits a set of generators $\{\gamma_1,\ldots,\gamma_n\}$ such that for all $r\geq 1$, $1\leq i_1<\cdots<i_r\leq n$, the traces $\mathrm{Tr}_g(\gamma_{i_1}\cdots\gamma_{i_r})$ are all geometric Salem numbers. Let $f_{i_1\cdots i_r}$ denote the minimal polynomial of $\mathrm{Tr}_g(\gamma_{i_1}\cdots\gamma_{i_r})$ over $\mathbb{Q}$. Then $(\Sigma,g)$ is marked length variety rigid with respect to the set of polynomials 
   $$\{X_1X_2-X_3-X_4\}\cup\{f_{i_1\cdots i_r}|r\geq 1,1\leq i_1<\cdots<i_r\leq n\}.$$
\end{thm}

A typical class of surfaces satisfying the hypotheses of Theorem~\ref{A} is given by arithmetic surfaces derived from quaternion algebras. At present, no non-arithmetic examples are known. Therefore, it is an interesting open question whether the conditions of Theorem~\ref{A} imply arithmeticity of the surface. 

Nevertheless, we can show a special case:
\begin{thm}\label{TH: salem to arithmtic}
    Let $S_{1,1}$ be a one-punctured torus, with standard generators $a$, $b$ of $\pi_1(S_{1,1})$, and let $g$ be a complete hyperbolic metric on $S_{1,1}$. If all three numbers $\mathrm{Tr}_g(a)$, $\mathrm{Tr}_g(b)$ and $\mathrm{Tr}_g(ab)$ are geometric Salem numbers, then $g$ is arithmetic.
\end{thm}

The paper is organized as follows: Section~\ref{Pre} recalls the necessary background on arithmetic Fuchsian groups. Section~\ref{Proofs} constructs non-arithmetic examples of marked length variety rigid surfaces and proves Theorems ~\ref{A} and ~\ref{TH: salem to arithmtic}. 

\section{Preliminaries}\label{Pre}
A standard reference for hyperbolic two- and three-manifolds, including arithmetic manifolds, from an algebraic perspective, is \cite{maclachlan2003arithmetic}. We briefly recall the necessary background for this paper.

We denote by $\SL(2,\mathbb{R})$
the group of real $2\times2$ matrices with determinant one and by $\PSL(2,\mathbb{R})$ the
quotient group $\SL(2,\mathbb{R}))/\{\pm I_2\}$ where $I_2$ is the $2\times 2$ identity matrix. Let $p:\SL(2,\mathbb{R})\rightarrow\PSL(2,\mathbb{R})$ be the projection. A discrete subgroup of $\PSL(2,\mathbb{R})$ is called \emph{Fuchsian}.  

A lattice of a locally compact, second countable topological group $G$ is a discrete subgroup $\Gamma$ such that $G/\Gamma$ has finite Haar measure. A lattice is called uniform if $G/\Gamma$ is compact and nonuniform otherwise. 

\medskip 

Arithmetic Fuchsian groups are obtained in the following way, for example, see \cite{maclachlan1987commensurability}: Let $k$ be an algebraic totally real number field
 so that the $\mathbb{Q}$-isomorphisms of $k$ into
$\mathbb{R}$ are $\phi_1$, $\phi_2$, $\ldots$, $\phi_n$ where we take $\phi_1=\mathrm{Id}$ and $\phi_i(k)\subset \mathbb{R}$ for $i=2,4,...,n$. Let $A$ be a quaternion algebra over $k$, which is ramified at only one real place, and thus there is an isomorphism
$$\rho: A\otimes_\mathbb{Q}\mathbb{R}=\SL(2,\mathbb{R})\oplus \bfH\oplus \bfH\oplus\cdots\oplus\bfH,$$
where $\bfH$ denotes Hamilton's quaternions. Denote $P$ the projection to the first factor.

Let $\mathcal{O}$ be an order in $A$, and $\mathcal{O}^1$  denote the group of elements of reduced norm 1. Then $P\rho(\mathcal{O}^1)$ is a lattice of $\SL(2,\mathbb{R})$ and $pP\rho(\mathcal{O}^1)$ is a lattice in $\PSL(2,\mathbb{R})$. The class of
\emph{arithmetic Fuchsian groups} is all Fuchsian lattices commensurable with such groups $pP\rho(\mathcal{O}^1)$. In addition, we say that a Fuchsian group is \emph{derived from a
quaternion algebra} if it is a subgroup of finite index in some $pP\rho(\mathcal{O}^1)$.

Among non-uniform arithmetic Fuchsian lattices, there is only one commensurability class, namely that of the modular group $\PSL(2, \mathbb Z)$. This is well known and follows, for instance, from \cite{maclachlan2003arithmetic}*{Theorem 8.1.2}.
Indeed, by \cite{maclachlan2003arithmetic}*{Theorem 8.1.2}, the associated quaternion algebra $A$ splits over the field $k$, and hence over every real embedding of $k$. It follows that $k$ has only one real embedding, so $k=\mathbb{Q}$, and consequently 
$A=\SL(2,\mathbb Q)$.

\subsection{Characterization of arithmetic Fuchsian and Kleinian groups}
Takeuchi characterizes arithmetic Fuchsian groups within the class of all Fuchsian lattices in \cite{takeuchi1975characterization}. 

Let $\Gamma$ be a Fuchsian group, and let $\Gamma^{(2)}$ denote the subgroup generated by the squares of elements of $\Gamma$. Note that if $\Gamma$ is finitely generated then $\Gamma^{(2)}$ is of finite index in $\Gamma$.
\begin{thm}[\cite{takeuchi1975characterization}, \cite{borel1981commensurability}]\label{arithmetic}
    If $\Gamma$ is an arithmetic Fuchsian group, then $\Gamma^{(2)}$ is derived from a quaternion algebra.
\end{thm}
\begin{thm}[\cite{takeuchi1975characterization}]\label{THM: Fuchsian}
Let $\Gamma$ be a Fuchsian lattice. Then $\Gamma$ is derived from a quaternion algebra over a totally real algebraic number field if and only if $\Gamma$ satisfies the following two conditions:
\begin{enumerate}
    \item 
$K:=\mathbb{Q}(\mathrm{Tr}(\Gamma))$ is an algebraic number field of finite degree and $\mathrm{Tr}(\Gamma)$
is contained in the ring of integers $\mathcal{O}_K$ of $K$.
\item
For any embedding $\phi$ of $K$ into $\mathbb{C}$, which is not the identity, $\phi(\mathrm{Tr}(\Gamma))$
is bounded in $\mathbb{C}$.
\end{enumerate}
\end{thm}

\section{Marked length variety and rigidity}\label{Proofs}
In this section, we construct non-arithmetic surfaces, both closed and non-compact, that are marked length variety rigid. We then prove Theorems~\ref{A} and \ref{TH: salem to arithmtic} after presenting these examples.
\subsection{A rigid, non-arithmetic one-punctured surface}

Consider a surface $S_{1,1}$ of genus one with one puncture, whose fundamental group is $\pi_1(S_{1,1})=F(a,b)$. We denote by $T_{1,1}$ the Teichm\"{u}ller space of $S_{1,1}$. A point $g\in T_{1,1}$ is uniquely determined by the Fricke coordinates $(x,y,z)=({\mathrm{Tr}}_g(a),{\mathrm{Tr}}_g(b),{\mathrm{Tr}}_g(ab))$ which satisfy the relation
\begin{equation}\label{ME}
x^2+y^2+z^2-xyz=0, \quad x>2, y>2, z>2.
\end{equation}

These are called Fricke coordinates, as in \cite{keen1977rough}. Alternative conventions, based on matrix entries, also appear in the literature \cite{keen1977rough}. 

Now, let $g\in T_{1,1}$ satisfy the marked length pattern:
$$\mathrm{Tr}_g(a)=\mathrm{Tr}_g(b),\quad \mathrm{Tr}_g(a^2)=\mathrm{Tr}_g(a^2b).$$
Since $\mathrm{Tr}_g(a^2)=x^2-2$ and $\mathrm{Tr}_g(a^2b)=zx-y$, these conditions are equivalent to
\begin{equation}\label{1}
   x=y, \quad x^2-2=zx-y. 
\end{equation}
Then $g$ defines a non-arithmetic surface with marked length variety rigidity, as shown by the following lemmas.

\begin{lemma}
The metric $g$ is unique.
\end{lemma}

\begin{proof} 
From equation~\eqref{1}, we have $z=\frac{x^2+x-2}{x}$. Substituting this into equation~\eqref{ME}, we obtain
\[
    x^2+x^2+\left(\frac{x^2+x-2}{x}\right)^2-x^2\left(\frac{x^2+x-2}{x}\right)=0.
\]
Since $x\neq 0$, multiplying through by $x^2$ and simplifying yields
\begin{equation}\label{eq1}
x^5-2x^4-4x^3+3x^2+4x-4=0.
\end{equation}
This equation has a unique real root $x_0\approx 2.91330$. For this value, both $y$ and $z$ are greater than 2, so all the required conditions are satisfied. Hence $g$ is the unique hyperbolic metric satisfying these conditions.
\end{proof}
\begin{lemma}
The metric $g$ is non-arithmetic.
\end{lemma}
\begin{proof}
Suppose, for contradiction, that $g$ is arithmetic. Then its fundamental group is commensurable with $\PSL_2(\mathbb Z)$. Consequently, there exist $k,m\in \mathbb N^+$ such that $\mathrm{Tr}_g(a^k)=m$. Let $\lambda=\left(\frac{m+\sqrt{m^2-4}}{2}\right)^{\frac{1}{k}}$, which is an eigenvalue of $a$. Then we have $x_0=\lambda+\frac{1}{\lambda}$. In particular, $x_0$ is solvable by radicals. 

However, the polynomial appearing on the left-hand side of equation~\eqref{eq1} is irreducible over $\mathbb{Q}$ and is the minimal polynomial of $x_0$ over $\mathbb Q$. Its Galois group is $S_5$, which is not solvable. This contradicts the fact that $x_0$ is solvable by radicals. Therefore, the metric $g$ is not arithmetic.
\end{proof}

\begin{lemma}\label{B}
The metric $g$ is marked length variety rigid.
\end{lemma}
\begin{proof}Let $S=\{X_1-X_2, X_1X_2-X_3-X_4\}$. We claim that $g$ is marked length variety rigid with respect to $S$.

Let $g'$ be a pinched, negatively curved, complete metric on $S$ such that $V_S^g\subset V_S^{g'}$.
Since 
$$\mathrm{Tr_g}(\gamma_1)\mathrm{Tr_g}(\gamma_2)-\mathrm{Tr_g}(\gamma_1\gamma_2)-\mathrm{Tr_g}(\gamma_1\gamma_2^{-1})=0$$
for all $\gamma_1,\gamma_2\in \bar{\pi}_1(S_{1,1})$, the same relation holds for $g'$. Following the algorithm outlined in \cite{horowitz1972characters}*{section 3}, this implies that $\mathrm{Tr}_{g'}$ is completely determined by $x_0=\mathrm{Tr}_{g'}(a)$, $y_0=\mathrm{Tr}_{g'}(b)$ and
$z_0=\mathrm{Tr}_{g'}(ab)$. Note that $[a,b]$ is a parabolic element, hence 
\[
\mathrm{Tr}_{g'}([a,b])=2.
\]
In particular,
$\mathrm{Tr}_{g'}(-1,[a,b])=-2$ is equivalent to
$$x_0^2+y_0^2+z_0^2-x_0y_0z_0=0.$$
Let $\bar{g}$ be the hyperbolic metric determined by $(x_0,y_0,z_0)$. Applying the same algorithm as before, we then have $\mathrm{Tr}_{g'}=\mathrm{Tr}_{\bar{g}}$, i.e. $\ell_{g'}=\ell_{\bar{g}}$. Hence, $g'=\bar{g}$ by the marked length spectrum rigidity.

Since $g$ is the unique point in $T_{1,1}$ satisfying all the required conditions, we conclude $g=\bar{g}=g'$.
\end{proof}

\subsection{A rigid, non-arithmetic genus-two surface} The construction and argument in this subsection are similar to those in the previous subsection.

We begin by recalling the structure of the Teichm\"{u}ller space $T_{2,0}$ of genus-two surfaces.
\begin{thm}[\cite{keen1977rough}*{Theorem 5.1}]\label{TH: T_{2,0}}
    The Teichm\"{u}ller space $T_{2,0}$ of genus two surfaces is described by the six-manifolds in $\mathbb{R}^9$ by the equations:
    \begin{equation}\label{eq: T20-1}
        x_1^2+x_2^2+x_3^2-x_1x_2x_3+x_7=2,
    \end{equation}
    \begin{equation}\label{eq:T20-2}
        x_4^2+x_5^2+x_6^2-x_4x_5x_6+x_7=2,
    \end{equation}
       \begin{equation}\label{eq: T20-3}
           x_7^2+x_8^2+x_9^2-x_7x_8x_9+(x_1^2+x_4^2)x_7+2x_1x_4(x_8+x_9)+2x_1^2+2x_4^2+x_1^2x_4^2=4,
       \end{equation}
    where $x_i$, $1\leq i
    \leq 9$ denote the traces of a set of closed geodesics, and $x_i>2$ for all $1\leq i\leq 9$.
\end{thm}
Let $S_{2,0}$ be a surface of genus two, and let $\gamma_i\in \pi_1(S_{2,0})$, $1\leq i
\leq 9$, be such that $\ell_g(\gamma_i)=x_i$ for all $g\in T_{2,0}$ as in \cite{keen1977rough}.

Now, let $g\in T_{2,0}$ satisfy the marked length pattern
\begin{equation}\label{eq: condition 1}
    x_1=x_2=x_3=x_4=x_5=x_6=x_8,
\end{equation}
\begin{equation}\label{eq: condition 2}
    x_9=\mathrm{Tr}(\gamma_1^2)=x_1^2
-2.
\end{equation}
We will show that $g$ is a non-arithmetic, marked length variety rigid metric. The proof proceeds in three lemmas.
\begin{lemma}
    The metric $g$ is unique.
\end{lemma}
\begin{proof}
   Let $x_1=x_2=x_3=x_4=x_5=x_6=x_8=x$. By equations~\eqref{eq: T20-1} and ~\eqref{eq: condition 1}, we have
   \[
   x_7=x^3-3x^2+2. 
   \]
   Since $x_9=x^2-2$, equation~\eqref{eq: T20-3} is equivalent to 
   \begin{equation}\label{eq: the equation for 20}
       x^5-9x^4+2x^3+11x^2-4x-4=0.
   \end{equation}
 Equation~\eqref{eq: the equation for 20} has three real solutions $x\approx 8.62734$, $x\approx -0.943951$, $x\approx -0.510887$ and two complex roots $x\approx -0.91375\pm 0.35563i.$ The metric $g$ corresponds to the solution $x\approx 8.62734.$ Since this is the only solution satisfying the required conditions, the metric $g$ is unique.
\end{proof}
\begin{lemma}
    The metric $g$ is non-arithmetic.
\end{lemma}
\begin{proof}
    Let $f(x)=x^5-9x^4+2x^3+11x^2-4x-4$. Since $f$ has no integer roots, it has no linear factor over $\mathbb{Q}$. Moreover, the sum of any two real roots of $f$ is not an integer, and the real part of each complex root is not a half-integer. Hence $f$ has no quadratic factor over $\mathbb{Q}$. Since $\mathrm{deg}f=5$, any nontrivial factorization over $\mathbb{Q}$ would contain a linear or quadratic factor. Therefore, $f(x)$ is irreducible over $\mathbb{Q}$.

For an arithmetic Fuchsian group derived from a quaternion algebra, the trace of every element lies in a totally real field (see Section~\ref{Pre}). By Theorem~\ref{arithmetic}, for an arithmetic metric $g'$, and for any $\gamma\in \Gamma$, the number $\mathrm{Tr}_{g'}(\gamma^2)$ lies in a totally real field. 

Let $\alpha\approx 8.62734$ be a root of $f(x)$ corresponding to the metric $g$. For this metric, we have
\[
\mathrm{Tr}_g(\gamma_1^2)=\alpha^2-2,
\]
which is not totally real. Hence, the metric $g$ is not arithmetic.
\end{proof}

Before proving that the metric $g$ is marked length variety rigid, we first introduce a family of polynomials.

We begin by fixing notation. The fundamental group $\Gamma:=\pi_1(\Sigma_{2,0})$ admits the standard representation:
    \[\langle A_1, B_1, A_2, B_2\mid B_2^{-1}A_2^{-1}B_2A_2B_1^{-1}A_1^{-1}B_1A_1=1\rangle.\]
    Then we have $x_i$, $1\leq i\leq 9$ in Theorem~\ref{TH: T_{2,0}}, denote the signed trace of \[
    A_1, B_1, A_1B_1, A_2, B_2, A_2B_2, B_2^{-1}A_2^{-1}B_2A_2, A_2B^{-1}_1A_1^{-1}B_1, A_2A_1,\] respectively (see \cite{keen1977rough}*{p. 1212}).

    We associate a polynomial to each element of the following set: 
    \[G_{2,3}=\{A_1B_1, A_1A_2, A_1B_2, B_1A_2, B_1B_2, A_2B_2, A_1B_1A_2, A_1B_1B_2, A_1A_2B_2, B_1A_1A_2\}.\]
    To prepare for the explicit construction of these polynomials, we first analyze the associated Fuchsian representation.
    
    Let $r\approx 8.62734$ be a solution of equation~\eqref{eq: the equation for 20}. Since \[\mathrm{Tr}_g(A_2)=\mathrm{Tr}_g(B_2)=\mathrm{Tr}_g(A_2B_2)=r,\] up to conjugation, the metric $g$ admits a Fuchsian representation $i:\pi_1(S_{2,0})\to\SL(2,\mathbb{R})$ such that
    \[
    i(A_2)=\begin{pmatrix}
        \lambda &0\\
          0 &\lambda^{-1}
        \end{pmatrix},\quad i(B_2)=\begin{pmatrix}
        \frac{r}{\lambda+1} &1\\
           \frac{\lambda r^2}{(\lambda+1)^2}-1&\frac{\lambda r}{\lambda+1}
        \end{pmatrix},
    \]
    where $\lambda>1$ satisfies $\lambda+\lambda^{-1}=r$. 
    
    There exists a representation $i'$, belonging to a different conjugacy class, with the same trace, such that
    \[
    i'(A_2)=\begin{pmatrix}
        \lambda &0\\
          0 &\lambda^{-1}
        \end{pmatrix},\quad i'(B_2)=\begin{pmatrix}
        \frac{r}{\lambda+1} &-1\\
           1-\frac{\lambda r^2}{(\lambda+1)^2}&\frac{\lambda r}{\lambda+1}
        \end{pmatrix},
    \]
    which corresponds to the same metric $g$ with the opposite orientation.

    Let $E=i(B_2^{-1}A_2^{-1}B_2A_2)$. By \cite{maclachlan2003arithmetic}*{equation (3.5), p.~121},
    \[\mathrm{Tr}(E)=3r^2-r^3-2<0.\] 
    Since $E$ is hyperbolic, there exist $\mu<-1$ and a matrix $P\in\SL(2,\mathbb{R})$ such that
    \[
    E=P\begin{pmatrix}
        \mu& 0\\
        0&\mu^{-1}
    \end{pmatrix}P^{-1}.
    \]
    Furthermore, the entries of $P$ can be chosen in the number field $\mathbb{Q}(r,\lambda,\mu)$.

    For a real number $s$, denote 
    \[
    E(s)=P\begin{pmatrix}
        0& (-\mu)^s\\
        -(-\mu)^{-s}& 0
    \end{pmatrix}P^{-1}.
    \]

Note that the subgroups generated by $A_1$, $B_1$, and by $A_2$, $B_2$ are both fundamental groups of a genus-one surface with one geodesic boundary component.
By \cite{keen1977rough}*{Theorem 3.1} and the assumption on the traces, there exists a matrix $Q\in \SL(2,\mathbb{R})$ such that 
\[
Qi(A_2)Q^{-1}=i(A_1), Qi(B_2)Q^{-1}=i(B_1).
\]
Here, \cite{keen1977rough}*{Theorem 3.1} applies only to representations into $\PSL(2,\mathbb{R})$. As noted earlier, there are two conjugacy classes of representations into $\SL(2,\mathbb{R})$ with the same traces. The existence of $Q$ follows from the fact that both ordered pairs $(A_1, B_1)$ and $(A_2, B_2)$ determine the same orientation on the surface.
    
     Since $A_2^{-1}B_2^{-1}A_2B_2=B_1^{-1}A_1^{-1}B_1A_1\in\Gamma$, we have \[QEQ^{-1}=i(B_1^{-1}A_1^{-1}B_1A_1)=i(A_2^{-1}B_2^{-1}A_2B_2)=E^{-1}.\] 
     In particular, $Q$ interchanges the two fixed points of $E$. Hence $Q$ plays the same role as $R(x)$ in \cite{keen1973correction}*{second paragraph, p. 61} for the metric $g$, up to a difference in notation.  Since $E$ is hyperbolic, $Q$ is in the family $\pm E(s)$. Hence, there exists $s_0\in\mathbb{R}$ such that
    \[
    i(A_1)=E(s_0)i(A_2)E(s_0)^{-1}; \quad i(B_1)=E(s_0)i(B_2)E(s_0)^{-1}.
    \]

    Let \[
    \phi(s)=i(A_2)E(s)i(B_2^{-1}A_2^{-1}B_2)E(s)^{-1}; \quad \varphi(s)=i(A_2)E(s)i(A_2)E(s)^{-1}.
    \]
    Note that both $\phi(s)$ and $\varphi(s)$ are in the form
    \[A\begin{pmatrix}
        0& (-\mu)^s\\
        -(-\mu)^{-s}& 0
    \end{pmatrix}B\begin{pmatrix}
        0& (-\mu)^s\\
        -(-\mu)^{-s}& 0
    \end{pmatrix}^{-1}C\]
    for some matrices $A$, $B$ and $C$. Note that 
    \[
    \begin{pmatrix}
        0& (-\mu)^s\\
        -(-\mu)^{-s}& 0
    \end{pmatrix}\begin{pmatrix}
        a_{11}&a_{12}\\
        a_{21}& a_{22}
    \end{pmatrix}\begin{pmatrix}
        0& (-\mu)^s\\
        -(-\mu)^{-s}& 0
    \end{pmatrix}^{-1}=\begin{pmatrix}
        a_{22}&-(-\mu)^{2s}a_{21}\\
        -(-\mu)^{-2s}a_{12}& a_{11}
    \end{pmatrix}.
    \]
    Then there exist $t_{ij}\in \mathbb{Q}(r,\lambda,\mu)$, for $1\leq i\leq 2$ ,$1\leq j\leq3$ such that
    \[
    \mathrm{Tr}(\phi(s))=t_{11}+t_{12}(-\mu)^{2s}+t_{13}(-\mu)^{-2s},
    \]
    \[
    \mathrm{Tr}(\varphi(s))=t_{21}+t_{22}(-\mu)^{2s}+t_{23}(-\mu)^{-2s}.
    \]
    
    By the definition of $s_0$, we have $\left|\mathrm{Tr}(\phi(s_0))\right|=x_8$ and $\left|\mathrm{Tr}(\varphi(s_0))\right|=x_9$. Therefore, 
    \[
    \mathrm{Tr}(\phi(s_0))=\pm r; \quad \mathrm{Tr}(\varphi(s_0))=\pm (r^2-2),
    \]
and $s_0$ is uniquely determined by these two conditions. It follows that $(-\mu)^{2s_0}$ has degree at most two over the field $\mathbb{Q}(r,\lambda,\mu)$. Note that the sign of $\mathrm{Tr}(\phi(s_0))$ and $\mathrm{Tr}(\varphi(s_0))$ is irrelevant for the argument. 

Let $X\in G_{2,3}$.
Similarly, there exist $t_{X,i}\in \mathbb{Q}(r,\lambda,\mu)$ such that
\[
\mathrm{Tr}(i(X))=t_{X,1}+t_{X,2}(-\mu)^{2s_0}+t_{X,3}(-\mu)^{-2s_0}.
\]
It follows that $\mathrm{Tr}(i(X))\in \mathbb{Q}(r,\lambda,\mu, (-\mu)^{2s_0}).$

Let $\sigma$ be a field embedding of $\mathbb{Q}(r,\lambda,\mu, (-\mu)^{2s_0})$ into $\mathbb{R}$ such that $\sigma(r)=r$. Since $\sigma$ fixes the traces of the nine elements corresponding to the
coordinates $x_i$ of the Teichm\"{u}ller space ($1\leq i\leq 9$), it follows that the representation $\sigma\circ  i$ has the same nine trace
coordinates as $i$. Therefore 
$\sigma\circ  i$ is conjugate to $i$. Consequently, for any element $X$, we have $\sigma(\mathrm{Tr}(i(X)))=\mathrm{Tr}(i(X))$. Hence $\mathrm{Tr}(i(X))\in \mathbb{Q}(r)$. 

For each $X$, let $f_X$ be the polynomial of degree less than 6 such that $\mathrm{Tr}(i(X))=f_X(r).$ Define 
\[h_X=X_1-f_X(X_2).\]
Now let $S'$ be the set of polynomials consisting of \[X_1-X_2,\quad X_1X_2-X_3X_4,\quad h_X\ \  (X\in G_{2,3}),\] together equation~\eqref{eq: T20-3}.
\begin{lemma}\label{MLVR}
    The metric $g$ is marked length variety rigid with respect to $S'$.
\end{lemma}
The proof of Lemma~\ref{MLVR} is similar to that of Lemma~\ref{B}, with two modifications. First, note that equations~\eqref{eq: T20-1} and $\eqref{eq:T20-2}$ for the coordinates are consequences of the relation $X_1X_2-X_3-X_4$, and hence they remain valid in the present setting. Second, the algorithm outlined in \cite{horowitz1972characters}*{section 3} should be replaced by \cite{maclachlan2003arithmetic}*{Lemma 3.5.3}. 

We leave the details to the reader.
\subsection{Proof of Theorem~\ref{A} and ~\ref{TH: salem to arithmtic}}

\begin{proof}[Proof of Theorem~\ref{A}]
Let $S=\{X_1X_2-X_3-X_4\}\cup\{f_{i_1\cdots i_r}|r\geq 1,1\leq i_1<\cdots<i_r\leq n\}$.
    Assume that $g'$ is a complete, pinched, negatively curved metric on $\Sigma$ such that $V_S^g\subset V_S^{g'}$. 

    Then, for each $f_{i_1\cdots i_r}$, we have $V_{f_{i_1\cdots i_r}}^g\subset V_{f_{i_1\cdots i_r}}^{g'}$. Hence $\mathrm{Tr}_{g'}(\gamma_{i_1}\cdots\gamma_{i_r})$ is a root of $f_{i_1\cdots i_r}$. By assumption, all Galois conjugates of $\mathrm{Tr}_{g}(\gamma_{i_1}\cdots\gamma_{i_r})$ have absolute value less than 2. 
    
    On the other hand, since $\gamma_{i_1}\cdots\gamma_{i_r}$ is not elliptic in the isometry group of $(\widetilde M,g')$,  we have $\mathrm{Tr}_{g'}(\gamma_{i_1}\cdots\gamma_{i_r})\geq 2$. Hence $\mathrm{Tr}_{g'}(\gamma_{i_1}\cdots\gamma_{i_r})=\mathrm{Tr}_{g}(\gamma_{i_1}\cdots\gamma_{i_r})$.

    As in the proof of Lemma~\ref{B}, for all $\eta_1$, $\eta_2\in \bar{\pi}_1(\Sigma)$, we have
    \[
    \mathrm{Tr}_{g'}(\eta_1)\mathrm{Tr}_{g'}(\eta_2)-\mathrm{Tr}_{g'}(\eta_1\eta_2)-\mathrm{Tr}_{g'}(\eta_1\eta_2^{-1})=0.
    \]
    The algorithm in \cite{horowitz1972characters}*{section 3} then implies $\mathrm{Tr}_{g'}=\mathrm{Tr}_{g}$. Therefore, $g'=g$ by the marked length spectrum rigidity of surfaces. 
    
    This completes the proof of Theorem~\ref{A}.
\end{proof}
\begin{proof}[Proof of Theorem~\ref{TH: salem to arithmtic}]
    Let $x=\mathrm{Tr}_g(a)$, $y=\mathrm{Tr}_g(b)$, $z=\mathrm{Tr}_g(ab)$. Then $T_{1,1}$ is described in $\mathbb{R}^3$ by 
    \begin{equation}\label{T11}
        x^2+y^2+z^2-xyz=0,\quad x>2, y>2, z>2.
    \end{equation}

    Assume that $x$, $y$, and $z$ are geometric Salem numbers, and let $K=\mathbb{Q}[x,y,z]$. If $\mathrm{Deg}_{\mathbb{Q}}(x)>1$, then there exists a field embedding $\sigma: K\to \mathbb{R}$ such that $|\sigma(x)|<2$. Applying $\sigma$ to equation~\eqref{T11}, we obtain
    \[
    \sigma(z)^2-[\sigma(x)\sigma(y)]\sigma(z)+[\sigma(x)^2+\sigma(y)^2]=0.
    \]
    The discriminant of this quadratic is
    \[
    \Delta=\sigma^2(x)\sigma^2(y)-4\sigma^2(x)-4\sigma^2(y)=(\sigma^2(x)-4)\sigma^2(y)-4\sigma^2(x)<0.
    \]
    So $\sigma(z)$ is not real. This contradicts the fact that $z$ is a geometric Salem number. Therefore, $\mathrm{Deg}_{\mathbb{Q}}(x)=1$. Similarly, $\mathrm{Deg}_{\mathbb{Q}}(y)=\mathrm{Deg}_{\mathbb{Q}}(z)=1$. That is, all three numbers are integers. 
    
    Applying the algorithm in \cite{horowitz1972characters}*{section 3}, the traces of all closed geodesics on $S_{1,1}$ are integers. By Theorem~\ref{THM: Fuchsian}, $g$ is arithmetic.
\end{proof}

\subsection{Acknowledgments} The author thanks the referee for a careful reading of the manuscript and for numerous valuable comments and suggestions that significantly improved the presentation of this paper.

\end{document}